\documentclass{article}
\usepackage[utf8]{inputenc}

 \usepackage{tabu}

\usepackage{amsmath}
\DeclareMathOperator{\sign}{sign}

\usepackage{algorithm2e}
\usepackage{xcolor}

\SetAlFnt{\footnotesize}

\SetAlCapSty{xAlCapSty}

\SetCommentSty{xCommentSty}

\SetNlSty{mynlfont}{}{} 

\LinesNumbered

\SetSideCommentRight

\DontPrintSemicolon

\RestyleAlgo{algoruled}

\usepackage[driverfallback=dvipdfm,
bookmarks=true,    
bookmarksopen=true,
bookmarksopenlevel=\maxdimen,
bookmarksdepth=10
]{hyperref}
\hypersetup{
pdfcenterwindow=true, 
pdfpagelabels=true,
pagebackref=true,            
hypertexnames=true,
plainpages=false,
unicode=false,     
pdftoolbar=true,                
    pdfmenubar=true,          
    pdffitwindow=false,         
    pdfstartview={FitH},        
    pdftitle={On Third-Order Accuracy of MUSCL Scheme},
    pdfauthor={Hiroaki Nishikawa},    
    pdfsubject={On Third-Order Accuracy of MUSCL Scheme},       
    pdfcreator={Hiroaki Nishikawa},   
    pdfproducer={Hiroaki Nishikawa}, 
    baseurl={http://www.hiroakinishikawa.com},
    pdfkeywords={}, 
    pdfnewwindow=true,     
    colorlinks=true,       
    linkcolor=black,       
    citecolor=black,       
    filecolor=black,       
    urlcolor=black,        
  breaklinks=true,
  hyperfigures=true,
  backref=true,
breaklinks=false, 
pdfpagelabels,      
pagebackref,        
hypertexnames=true, 
plainpages=false,   
naturalnames        
}
\usepackage[all]{hypcap}

\usepackage[thinlines]{easytable}

\usepackage{color}      
\usepackage{placeins}
\usepackage{indentfirst}
\usepackage{multirow}

\usepackage{amssymb}
\usepackage{amsmath,latexsym}
\usepackage{url}
\usepackage{enumerate}
\usepackage{graphicx}

\usepackage{booktabs}

\usepackage{subcaption}
\usepackage{cleveref}
\usepackage{mathtools}

\usepackage{blindtext,graphicx}
\usepackage[absolute]{textpos}
\begin{document}


\title{\bf On Estimating Machine-Zero Residual}

\author{ 
{Hiroaki Nishikawa}\thanks{ {hiro@nianet.org}  }\\
  {\normalsize\itshape National Institute of Aerospace, Hampton, VA 23666, USA}
}

\date{\today}
\maketitle

\begin{abstract}  
In this paper, we propose two techniques to estimate the magnitude of a machine-zero residual for a given problem, 
which is the smallest possible residual that can be achieved when we solve a system of discretized equations. 
We estimate the magnitude of the machine-zero residual by a norm of residuals computed 
with a randomly-perturbed approximate solution that is considered as close in magnitude to an exactly-converged solution. 
One method uses free-stream values as the approximate solution, and the other
uses a current solution during an iterative solve as the approximate solution via the method of manufactured solutions. 
Numerical results show that these estimates predict the levels of machine-zero residuals very accurately for all equations of the Euler and 
Navier-Stokes equations in a transonic flow over an airfoil and viscous flows over a cylinder and a flat plate.
\end{abstract}


\section{Introduction}

A machine-zero residual is defined as a residual (e.g., a finite-volume discretization of the Euler equations) with a converged solution 
substituted having machine-zero-level perturbation. It is the smallest possible residual value that a solver can achieve; it is 
 difficult to predict especially before a calculation. Also, it varies widely depending on the grid unit and the solution units especially 
when solving a dimensional form of flow equations. However, if it were possible to estimate a machine-zero residual before or during a calculation, then
it would be very useful. For example, if an initial solution happens to be a converged solution 
(e.g., a restart) or close to it, then a machine-zero residual estimate would help us determine that the residual is well satisfied already and thus 
we would stop even before calling a solver. On the other hand, if there is no estimate available, we will find that a solver generates very small solution updates 
but cannot tell if the solution is already converged (it could be stalling). One can look at the
corresponding residual level, but again one cannot tell whether it is the smallest possible value (it could be stalling).
Note again that the machine-zero residual level depends on many different factors and thus it can be very small or very large: 
e.g., being $1.0$e-09 does not necessarily mean a machine-zero residual. 

In practical fluid-dynamics solvers, iterative convergence is checked based on the reduction of the residual or iterative-solution-difference from their 
initial values. This is a practical thing to do, but it can be inefficient. For example, if an initial solution is very close to a converged solution (e.g., 
a solution at a previous time step in an inner iteration of an implicit time-stepping scheme with a very small time step), then the residual is already very 
small at the beginning of an iteration. Reducing the residual by two or three orders of magnitude may be much more than necessary; one 
order of magnitude reduction may suffice. If an accurate estimate of a machine-zero residual is available, then we can determine how well the discrete
equations are satisfied (e.g., five orders of magnitude above the machine-zero level) and stop the iteration earlier and could save computing time by a 
considerable amount over an entire unsteady simulation.

To the best of the author's knowledge at the time of writing this, no techniques are currently available for estimating a machine-zero residual level for 
a given flow problem. It is the objective of this paper to propose two such techniques and investigate how well they estimate machine-zero residual levels
for inviscid and viscous flow problems in two dimensions.

\section{Machine-Zero Residual}
Consider a nonlinear system of equations, which arise from a finite-volume-type discretization of the dimensional Euler/Navier-Stokes equations,
\begin{eqnarray}
  {\bf R}  ( {\bf U} ) = 0, \label{residual}
\end{eqnarray}
where $  {\bf R} $ and $  {\bf U} $ denote global vectors of residuals and numerical solutions in a given grid, 
which we solve by an iterative solver:
\begin{eqnarray}
 {\bf U}^{k+1} = {\bf U}^{k} + \Delta {\bf U}^{k},
 \end{eqnarray}
where $k$ is the iteration counter and $\Delta {\bf U}^{k}$ is a correction 
suggested by an iterative method (e.g., Newton's method). 
Let $ \overline{\bf U} $ be an exactly converged solution for a given problem with appropriate boundary conditions: 
\begin{eqnarray}
  {\bf R}  ( \overline{\bf U} ) = 0. \label{zero_residual}
\end{eqnarray}
In reality, the iteration will not find $\overline{\bf U}$ exactly but converge to $\overline{\bf U}$ with perturbation of the order of 
machine zero, $ \overline{\bf U} + \epsilon \overline{\bf U} $, where $\epsilon$ is a machine zero 
for a value of $O(1)$ (e.g., 1.0E-16 in double precision). This is the best we can hope for. 
Then, the minimum residual that can be attained, which we call the machine-zero residual, would be 
\begin{eqnarray}
   {\bf R}_{min}  =  {\bf R}  ( \overline{\bf U} + \epsilon \overline{\bf U} ),
\end{eqnarray}
which may be expanded as
\begin{eqnarray}
   {\bf R}_{min} 
   =  {\bf R}  ( \overline{\bf U} ) + \overline{ \frac{\partial {\bf R}}{\partial {\bf U}} } \epsilon \overline{\bf U} + O(\epsilon^2)
   =  \overline{  \frac{\partial {\bf R}}{\partial {\bf U}} } \epsilon \overline{\bf U} + O(\epsilon^2),
   \label{machine_zero_res}
\end{eqnarray}
where $\overline{ \frac{\partial {\bf R}}{\partial {\bf U}} } $ is the residual Jacobian evaluated at 
$ \overline{\bf U}$. It shows that the machine-zero residual depends on the magnitude of the solution as well as the Jacobian. Hence, 
there are various factors that affect the magnitude of the machine-zero residual: e.g., grid units and physical units.
Therefore, it is not necessarily of the order of 1.0E-16.

At this point, it is clear that the machine-zero residual cannot be computed because the exactly-converged solution $ \overline{\bf U}$ is not known. 

\section{Estimates of Machine-Zero Residual}

To estimate a machine-zero residual, we need to approximate the exactly-converged solution $\overline{\bf U}$. 
As our interest is only in estimating the magnitude of the machine-zero residual in some norm, the approximation just has to be
close in magnitude to the exactly-converged solution $\overline{\bf U}$. There are two possibilities. 
 
\subsection{Free-stream estimate $R_c$}  

One is a constant free-stream solution:
\begin{eqnarray}
  {R}_c
   =     || {\bf R}  ( {\bf U}_\infty +  \epsilon {r} {\bf U}_\infty ) ||, 
   \label{machine_zero_res_estimate_infty}
\end{eqnarray}
where ${r}$ is a random number ($0 \le  r \le 1$) and $ {\bf U}_\infty$ denotes a constant solution, which is an exact solution with a free-stream condition
applied at all boundaries: $ {\bf R}  ( {\bf U}_\infty ) = 0$. Note that in practice we apply different random numbers to solutions at different nodes 
(or cells) as we will show below. It is reasonable to set a free-stream condition at all boundaries, but we did not do so  
for the test cases shown in this paper.  
 This estimate can be computed before a calculation and also separately for different equations in a target system.

In our dimensional Euler/Navier-Stokes solver, we store the primitive variables ${\bf w} = ( p', u, v, T  ) $ at nodes, where $u$ and $v$ are velocity components, $T$ is the temperature, and $p'$ is the gauge pressure $p' = p- p_\infty$, where $p_\infty$ is a free-stream pressure. Hence, the residual is a function of the primitive variables stored at nodes: 
 $ {\bf R}  ( {\bf U}) =  {\bf R}  ( [ {\bf w}_1, {\bf w}_2, {\bf w}_3, \cdots,  {\bf w}_N ])$, 
 where $N$ denotes the number of nodes in a grid.  We define the perturbed free-stream solution at a node $j$ as
\begin{eqnarray}
  {\bf w}_\infty +  \epsilon r_j  {\bf w}_\infty
 = 
  \left [
\begin{array}{c}
\epsilon ( 1  +   \epsilon r_j  )
\\[1ex]
\mbox{maxmod}\left( u_\infty ,  \sign(u_\infty) \epsilon \right)  ( 1  +   \epsilon  r_j  ) 
\\[1ex]
\mbox{maxmod}\left( v_\infty , \sign(v_\infty)   \epsilon \right)  ( 1  +   \epsilon r_j  )
\\[1ex]
T_\infty ( 1  +   \epsilon r_j   )
\end{array}
\right ].
\end{eqnarray}
where $ r_j$ is a random number defined at the node $j$ and the maxmod function is defined by
\begin{eqnarray}
\mbox{maxmod}( a, b) =   \left \{
\begin{array}{rl}
a  , & \mbox{if} \,\, |a| \ge |b|, \\[3.5ex]
b , & \mbox{otherwise},
\end{array}
\right.
\end{eqnarray} 
Note that the free-stream value of the gauge pressure $p'$ is zero. 
The velocity components $u_\infty$ and $v_\infty$ are computed for a given set of a free-stream Mach number and an angle of attack.
It is emphasized again that the perturbed free-stream solution is different at different nodes. 
 
\subsection{Current solution as a manufactured solution $R_m$ }

The other possibility is to use a current solution ${\bf U}^k$ as an exact solution, which is possible by the method of manufactured solutions.
Let ${\bf S}$ be a source term defined by
\begin{eqnarray}
 {\bf S}= {\bf R}  ( {\bf U}^k ),
 \end{eqnarray}
Then, ${\bf U}^k$ is an exact solution to the problem,
\begin{eqnarray}
  {\bf R}  ( {\bf U} ) = {\bf S},
 \end{eqnarray}
which leads to the following estimate: 
\begin{eqnarray}
{ R}_m   =   ||  {\bf R}  ( {\bf U}^k +    \epsilon {r} {\bf U}^k) -{\bf S} ||.
\label{Rm}
\end{eqnarray}
This is a machine-zero residual for a slightly different problem, 
but expected to serve well as an estimate, at least in magnitude, for the original problem. In fact, if we expand it as 
\begin{eqnarray}
{ R}_m  \approx  \left| \left|    \frac{\partial {\bf R}}{\partial {\bf U}}   \epsilon {r} {\bf U}^k   \right| \right|, 
\end{eqnarray}
we find that the estimate $R_m$ is an approximation to the machine-zero residual in the form (\ref{machine_zero_res}). 
It would be an accurate estimate if the current solution is close to a converged solution in magnitude.
In practice, we again apply different random numbers to solutions at different nodes.  
As in the previous case, this estimate can also be computed for each equation in a target system. 
It may be updated at every iteration if needed. 
Note that this technique is more general than the free-stream estimate $R_c$; it can be used with a free-stream solution by
setting $ {\bf U}^k = {\bf U}_\infty$ (then there is no need to specify a free-stream condition at all boundaries).

In our dimensional solver, we define the perturbed solution in the primitive variables, at a node $j$, as
\begin{eqnarray}
 {\bf w}_j^k +  \epsilon  {r} {\bf w}_j^k 
 = 
  \left [
\begin{array}{c}
\mbox{maxmod}(    p'_j   ( 1  +   \epsilon    r_j    )  ,    \epsilon  r_j    ) 
\\[1ex]
\mbox{maxmod}(    u_j   ( 1  +   \epsilon    r_j    )  ,    \epsilon   r_j  ) 
\\[1ex]
\mbox{maxmod}(    v_j   ( 1  +   \epsilon   r_j    )  ,    \epsilon    r_j     ) 
\\[1ex]
T_j ( 1  +   \epsilon   r_j     )
\end{array}
\right ],
\end{eqnarray}
where the variables $(p'_j, u_j,v_j, T_j)$ are values at the $k$-th iteration.


 
\section{Results}
\label{sec:results}

In this section, we will present numerical results for inviscid and viscous test cases. The Euler and Navier-Stokes equations
are solved in the dimensional form. The discretization method is the edge-based method \cite{liu_nishikawa_aiaa2016-3969} with the Roe inviscid flux \cite{Roe_JCP_1981} and the alpha-damping viscous flux \cite{nishikawa:AIAA2010}: the residual is defined at a node $j$ by 
\begin{eqnarray}
{\bf Res}_j = \sum_{k \in \{ k_j\} } {\Phi}_{jk} ( {\bf n}_{jk} ) +  s_j V_j, 
\label{threed_fv_semidiscrete_system_00}
\end{eqnarray}
where $V_j$ is the measure of a dual control volume around the node $j$, $\{ k_j\}$ is a set of neighbor nodes 
of the node $j$, ${\Phi}_{jk}$ is the sum of inviscid and viscous numerical fluxes, and ${\bf n}_{jk}$ is the directed area vector
along the edge. The system of residual equations are then solved by an implicit defect-correction solver.  See Ref.\cite{liu_nishikawa_aiaa2016-3969} for further details. 
We store the primitive variables at nodes, reconstruct them to achieve second-order 
accuracy, and directly update them in a solver. 
For all problems, the free-stream pressure and temperature are defined by
\begin{eqnarray}
\quad p_\infty = 101325.0 \, \mbox{  [Pa] }, 
\quad
T_\infty = 288.15 \, \mbox{  [K] }, 
\end{eqnarray}
and the velocity will be determined for a given Mach number and an angle of attack.

 The machine-zero residual estimate $R_c$ is computed only once before
a calculation and $R_m$ is computed at every iteration to see how it changes. The estimates and actual residual norms will be computed in the $L_1$ norm 
for all cases (i.e., the arithmetic average of the absolute value of a residual over all nodes in a grid). The machine zero $\epsilon$ is set to be 
$\epsilon = 1.0$e-16 unless otherwise stated. 
We will show also the $L_1$ norm of the iterative solution error $dw(i)$ defined for each primitive variable $w_i$, $i=1,2,3,4$, at the $k$-th nonlinear iteration by
\begin{eqnarray}
 d w(i) = \frac{1}{N} \sum_{j \in \{nodes\} }  \frac{ |  w_j(i)^{k} - w_j(i)^{k-1}  |   }{  \tilde{w}(i)  } ,
\end{eqnarray}
where 
\begin{eqnarray}
 \tilde{w}(i) =   \left \{
\begin{array}{rl}
\displaystyle  \max_{j \in \{nodes\} } | w_j(i) |     ,  & \mbox{if $\displaystyle  \max_{j \in \{nodes\} } | w_j(i)  |   \ge $ 1.0e-05} , \\[3.5ex]
 1 , & \mbox{otherwise}.
\end{array}
\right.
\label{conventional_limiter}
\end{eqnarray} 
This quantity is computed here only to confirm machine-zero convergence of the solution variables. We would not determine
convergence based on it since a small iterative solution error does not necessarily mean convergence (e.g., it may be stalling). 
We wish to determine convergence by the residual norm, which tells us how well the discrete equations are satisfied.

\subsection{Inviscid transonic flow over a Joukowsky airfoil: $M_\infty=0.85$ at an angle of attack $1.25^\circ$}
\label{sec:results_airfoil} 
   
\begin{figure}[htbp!]
    \centering
      \begin{subfigure}[t]{0.48\textwidth}
        \includegraphics[width=\textwidth,trim=2 2 2 2,clip]{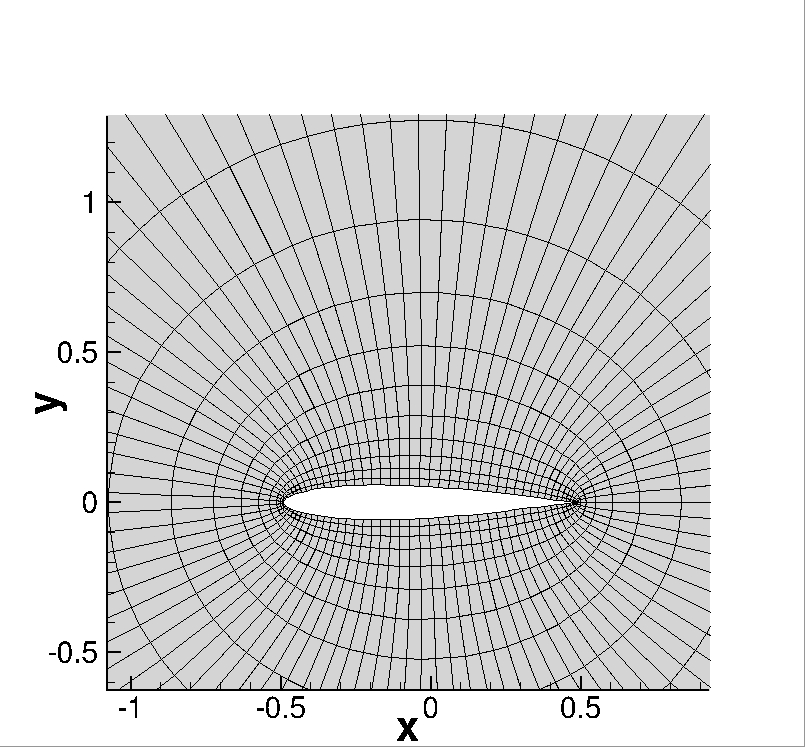}
          \caption{Grid.}
          \label{fig:inv_airfoil_grid}
      \end{subfigure}
      \hfill
      \begin{subfigure}[t]{0.48\textwidth}
        \includegraphics[width=\textwidth,trim=2 2 2 2,clip]{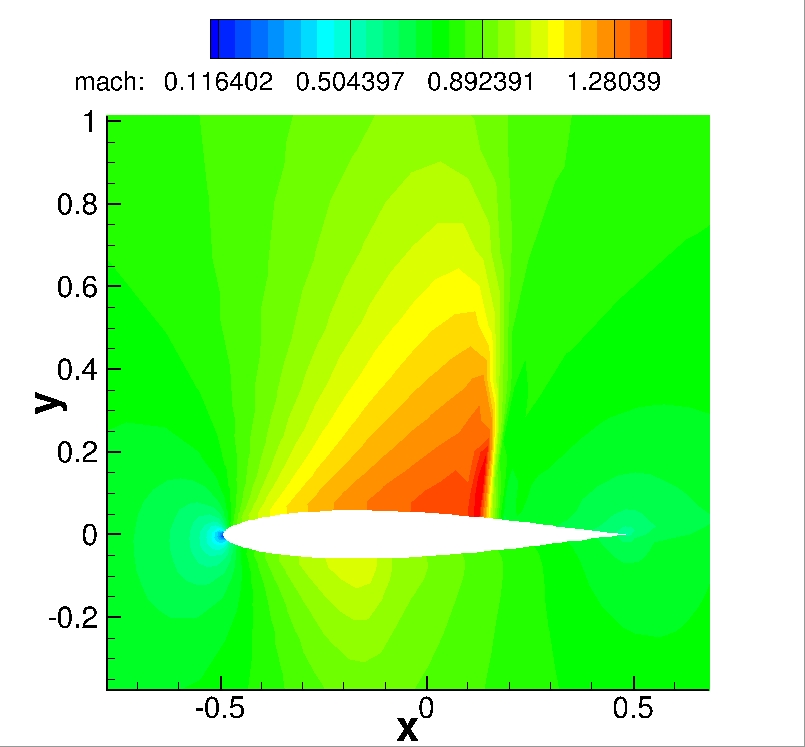}
          \caption{Mach contours.}
          \label{fig:inv_airfoil_press}
      \end{subfigure}
      \hfill
      \begin{subfigure}[t]{0.48\textwidth}
        \includegraphics[width=\textwidth,trim=2 2 2 2,clip]{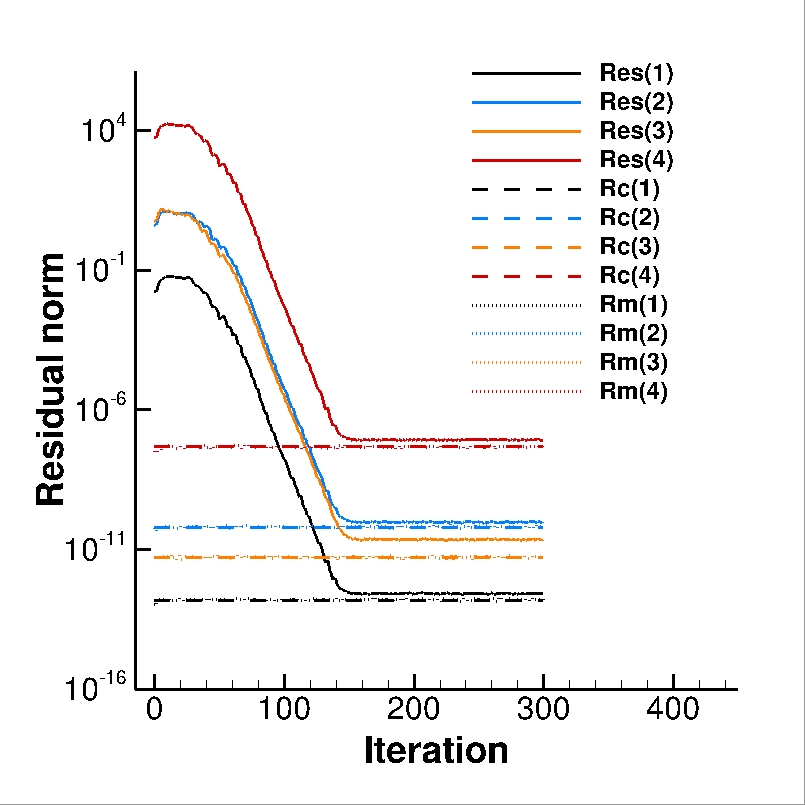}
          \caption{$L_1$ residual norm convergence.}
          \label{fig:inv_airfoil_res1}
      \end{subfigure}
      \hfill
      \begin{subfigure}[t]{0.48\textwidth}
        \includegraphics[width=\textwidth,trim=2 2 2 2,clip]{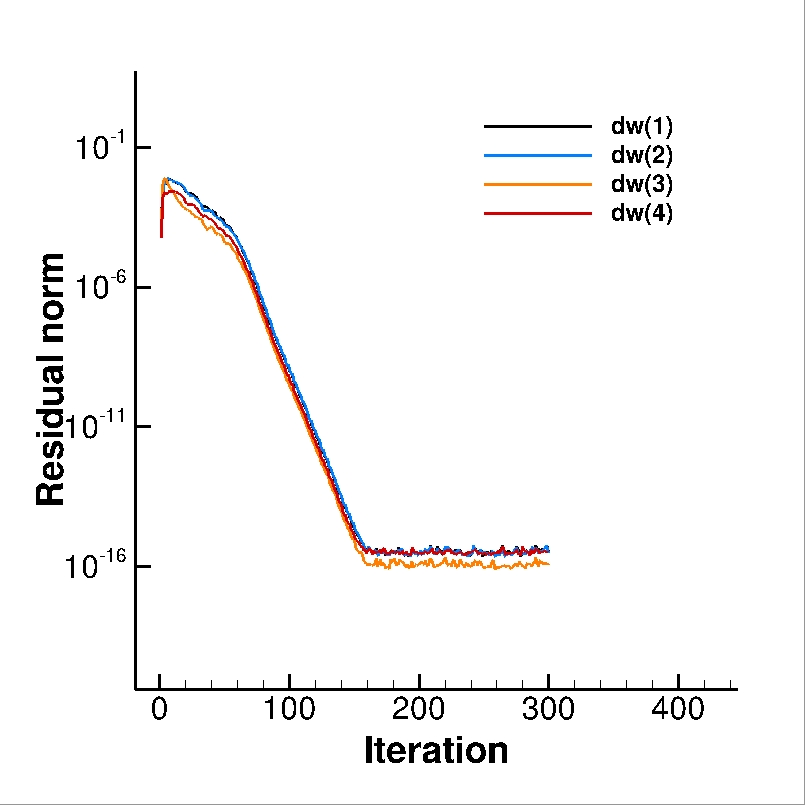}
          \caption{Iterative solution difference convergence.}
          \label{fig:inv_airfoil_dw}
      \end{subfigure}
      \caption{A transonic airfoil case. $Rc$ and $Rm$ computed for the four equations (1,2,3,4 = continuity, $x$-momentum, $y$-momentum, energy). 
      $R_m$ was updated at every iteration.
The actual residual norms, $Res(1), Res(2), Res(3)$, and $Res(4)$, stop going down nearly at the predicted levels. }
\label{fig:inv_airfoil} 
\end{figure}

We solve the Euler equations for a transonic flow over a Joukowsky airfoil at an angle of attack $1.25^\circ$ with $M_\infty=0.85$. 
Results are shown in Figure \ref{fig:inv_airfoil}.  
As can be seen, the machine-zero residuals are accurately predicted by both estimates. It is seen also that the estimate $R_m$ does not
change very much during the iteration.  
The iterative solution difference is also plotted. The results confirm that the solution is indeed converged to machine zero. 

\subsection{Viscous flow over a cylinder: $M_\infty=0.001$ and $Re_\infty=20$}
\label{sec:results_visc_cyl} 
   
\begin{figure}[htbp!]
    \centering
      \begin{subfigure}[t]{0.48\textwidth}
        \includegraphics[width=\textwidth,trim=2 2 2 2,clip]{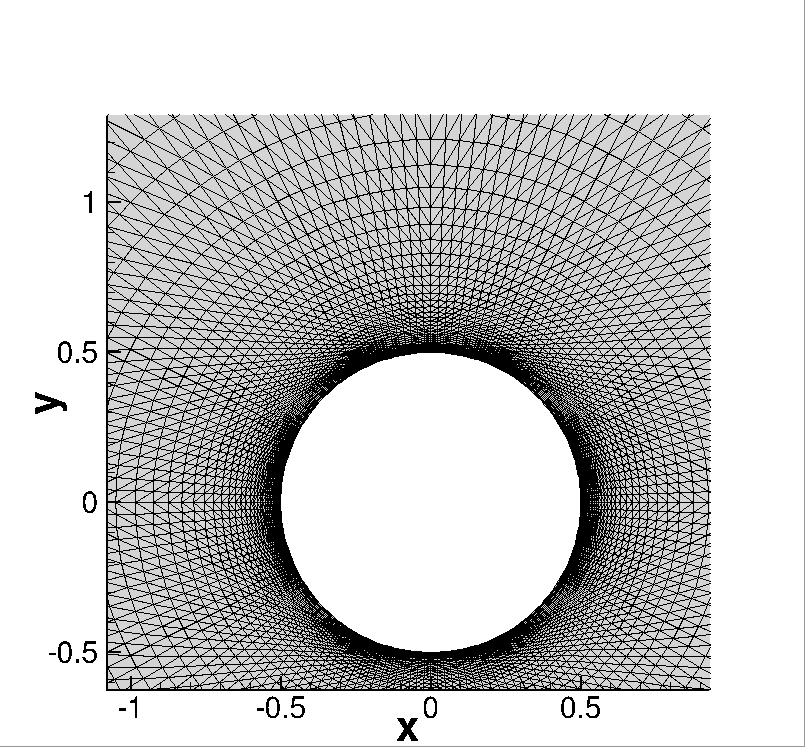}
          \caption{Grid.}
          \label{fig:inv_visc_cyl_grid}
      \end{subfigure}
      \hfill
      \begin{subfigure}[t]{0.48\textwidth}
        \includegraphics[width=\textwidth,trim=2 2 2 2,clip]{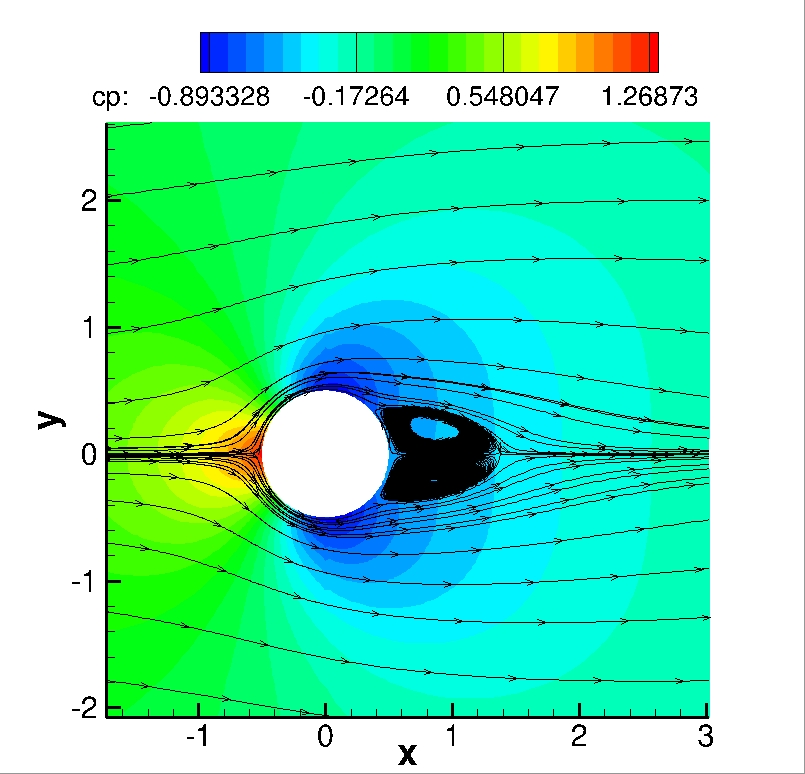}
          \caption{Pressure contours.}
          \label{fig:inv_visc_cyl_press}
      \end{subfigure}
      \hfill
      \begin{subfigure}[t]{0.48\textwidth}
        \includegraphics[width=\textwidth,trim=2 2 2 2,clip]{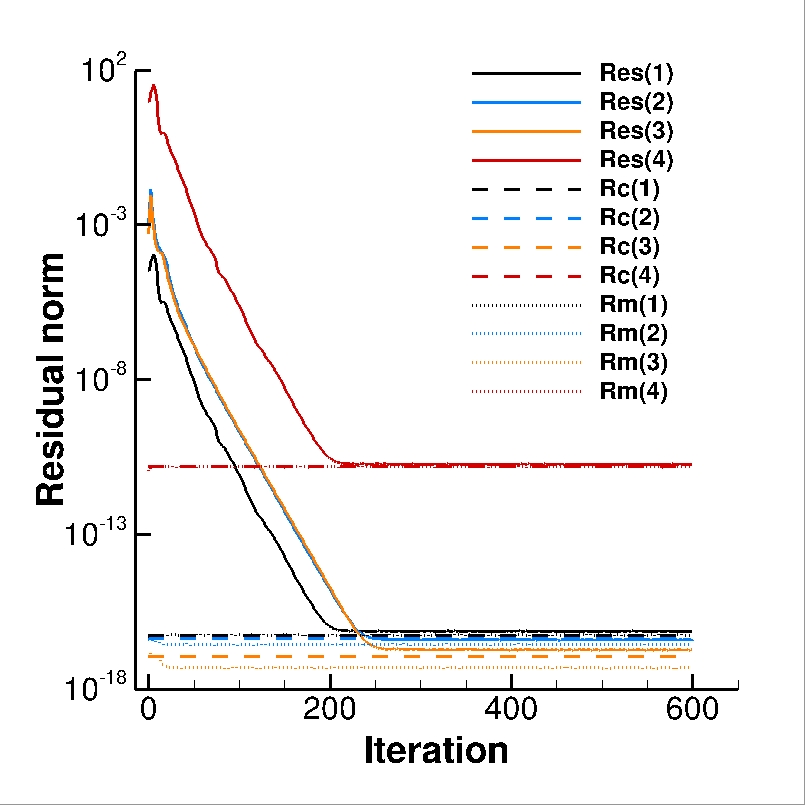}
          \caption{$L_1$ residual norm convergence.}
          \label{fig:inv_visc_cyl_res1}
      \end{subfigure}
      \hfill
      \begin{subfigure}[t]{0.48\textwidth}
        \includegraphics[width=\textwidth,trim=2 2 2 2,clip]{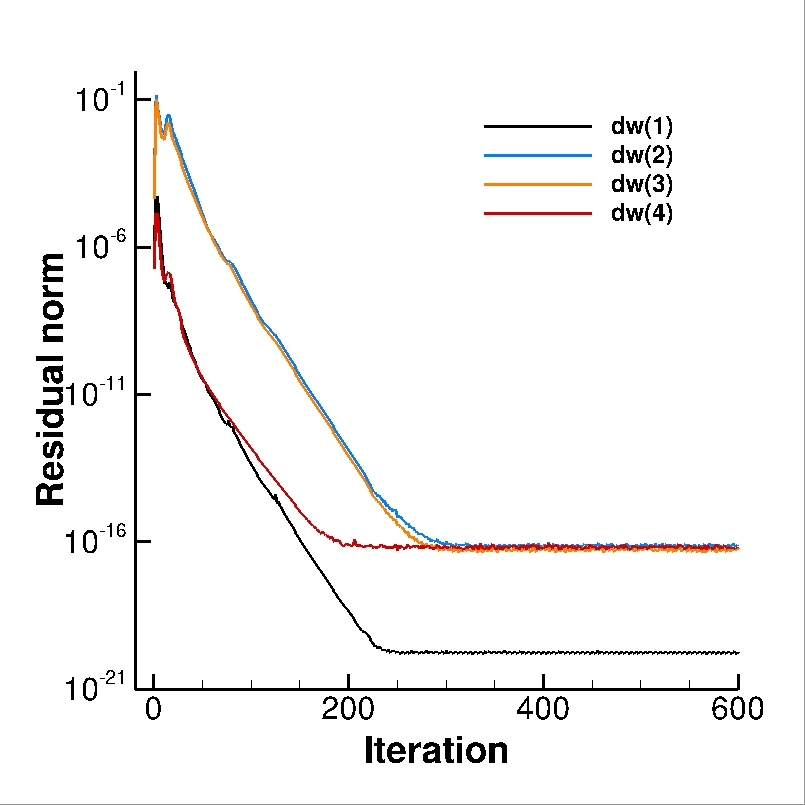}
          \caption{Iterative solution difference convergence.}
          \label{fig:inv_visc_cyl_dw}
      \end{subfigure}
      \caption{A viscous cylinder case.  $Rc$ and $Rm$ computed for the four equations (1,2,3,4 = continuity, $x$-momentum, $y$-momentum, energy). 
      $R_m$ was updated at every iteration.
The actual residual norms, $Res(1), Res(2), Res(3)$, and $Res(4)$, stop going down nearly at the predicted levels.}
\label{fig:inv_visc_cyl} 
\end{figure}

This is a low-Mach viscous flow test case. 
The free-stream velocity is determined for the Mach number 0.001 and a zero angle of attack. 
The viscosity is a constant and determined from $Re_\infty=20$. To overcome well-known low-Mach problems, we employed
here a preconditioned technique of Weiss and Smith \cite{WeissMaruszewskiSmith_AIAAJ1999} and solved an preconditioned 
Navier-Stokes system. Results are shown in Figure \ref{fig:inv_visc_cyl}. As can be seen, the estimates are very accurately predicted and 
the iterative solution difference has reached machine-zero convergence.

\subsection{Viscous flow over a flat plate: $M_\infty=0.15$ and $Re_\infty=10^4$}
\label{sec:results_visc_cyl} 

\begin{figure}[htbp!]
    \centering
      \begin{subfigure}[t]{0.48\textwidth}
        \includegraphics[width=\textwidth,trim=2 2 2 2,clip]{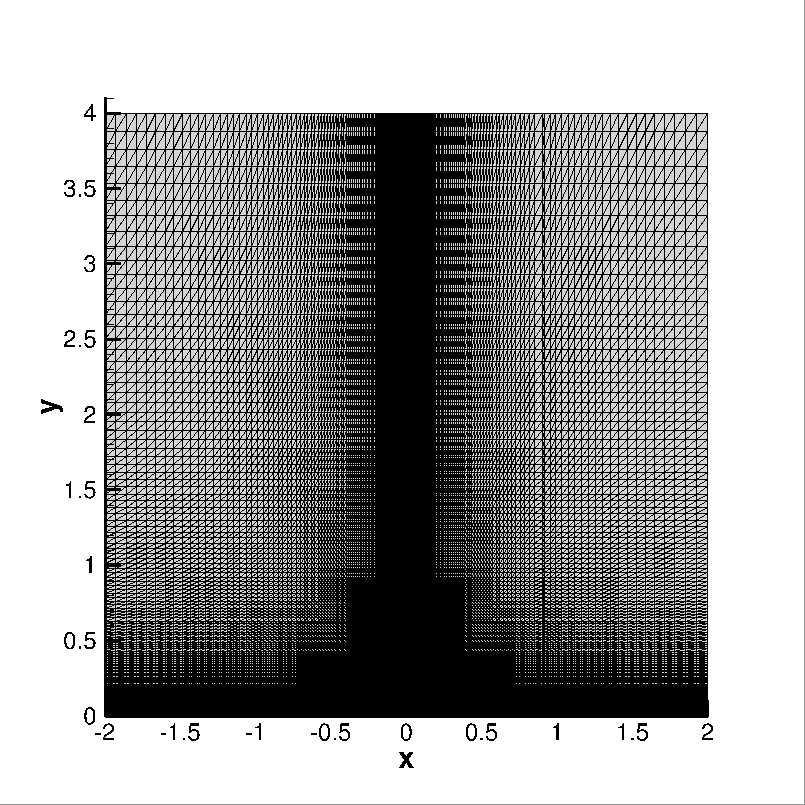}
          \caption{Grid (274$\times$194).}
          \label{fig:inv_visc_fp_grid}
      \end{subfigure}
      \hfill
      \begin{subfigure}[t]{0.48\textwidth}
        \includegraphics[width=\textwidth,trim=2 2 2 2,clip]{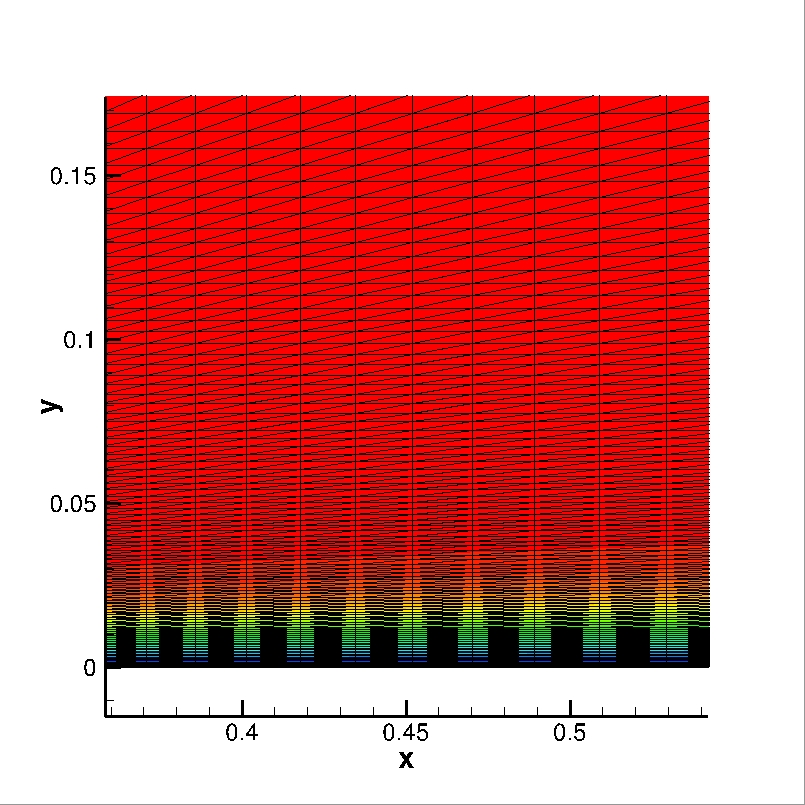}
                  \caption{Mach contours.}
          \label{fig:inv_visc_fp_press}
      \end{subfigure}
      \hfill
      \begin{subfigure}[t]{0.48\textwidth}
        \includegraphics[width=\textwidth,trim=2 2 2 2,clip]{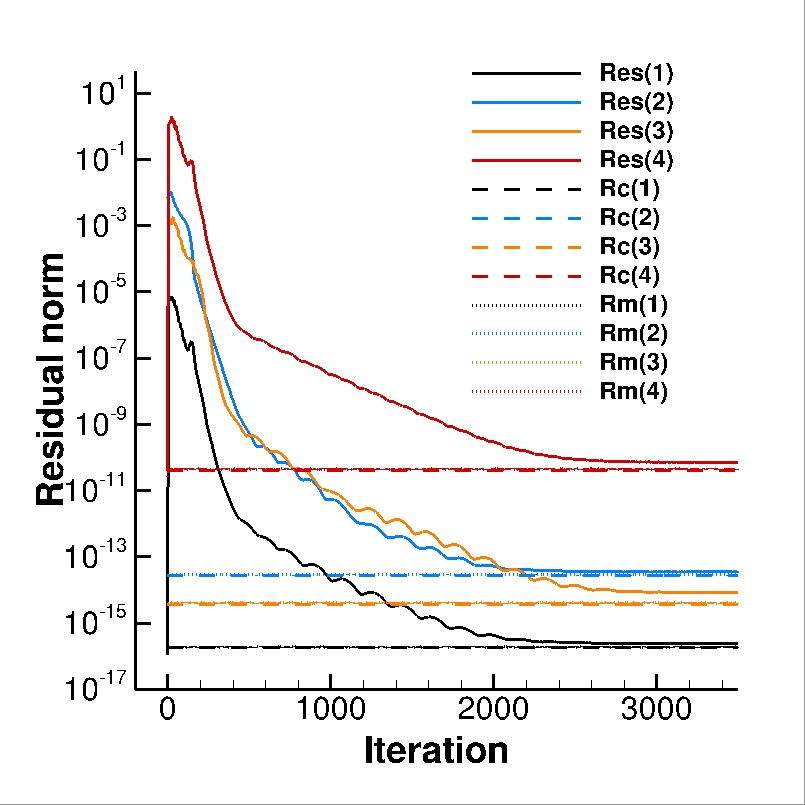}
          \caption{$L_1$ residual norm convergence.}
          \label{fig:inv_visc_fp_res1}
      \end{subfigure}
      \hfill
      \begin{subfigure}[t]{0.48\textwidth}
        \includegraphics[width=\textwidth,trim=2 2 2 2,clip]{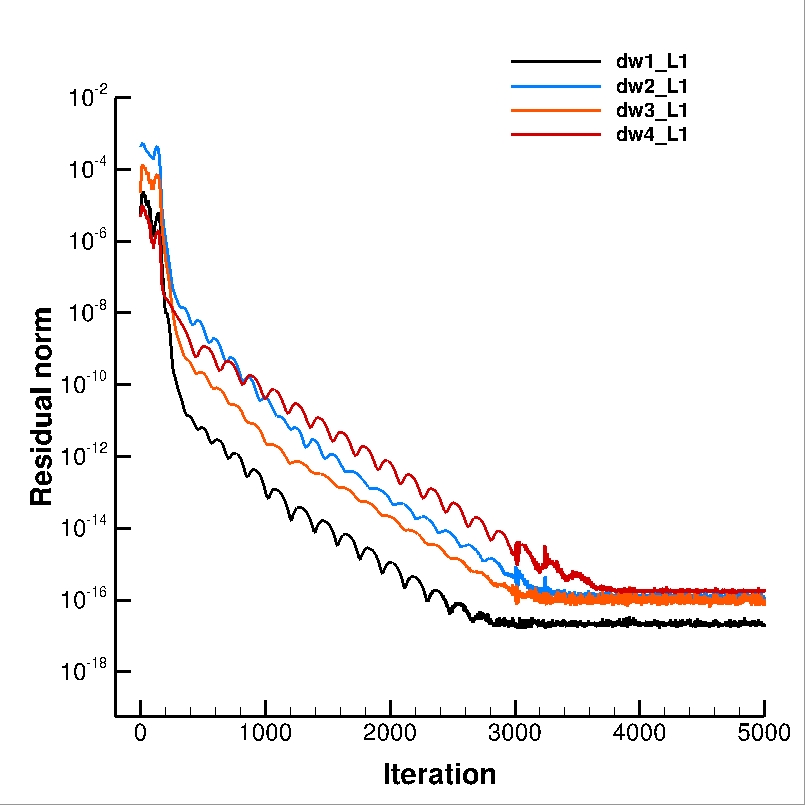}
          \caption{Iterative solution difference convergence.}
          \label{fig:inv_visc_fp_dw}
      \end{subfigure}
      \caption{A flat plate case.  $Rc$ and $Rm$ computed for the four equations (1,2,3,4 = continuity, $x$-momentum, $y$-momentum, energy). 
      $R_m$ was updated at every iteration.
The actual residual norms, $Res(1), Res(2), Res(3)$, and $Res(4)$, stop going down nearly at the predicted levels.}
\label{fig:inv_visc_fp} 
\end{figure}

We consider a viscous flow over a flat plate. We again solve the preconditioned Navier-Stokes equations in the dimensional form. 
The velocity and the viscosity are determined from $M_\infty=0.15$ (zero angle of attack) and $Re_\infty=10^4$. 
The grid is a triangular grid with 274$\times$194 nodes. 
Results are shown in Figure \ref{fig:inv_visc_fp}. It can be seen that the estimates very accurately predict the4 machine-zero residuals 
and the iterative solution difference has reached machine-zero convergence. 

\begin{figure}[htbp!]
    \centering
      \begin{subfigure}[t]{0.32\textwidth}
        \includegraphics[width=\textwidth,trim=2 2 2 2,clip]{./developer_test_res_estimate3/case_visc_fp_tria274x194_4eqns/tol_dw_10m16/myfig_res1_vs_itr.jpeg}
          \caption{$\epsilon=10^{-16}$.}
          \label{fig:inv_visc_fp_res1_16}
      \end{subfigure}
      \hfill
      \begin{subfigure}[t]{0.32\textwidth}
        \includegraphics[width=\textwidth,trim=2 2 2 2,clip]{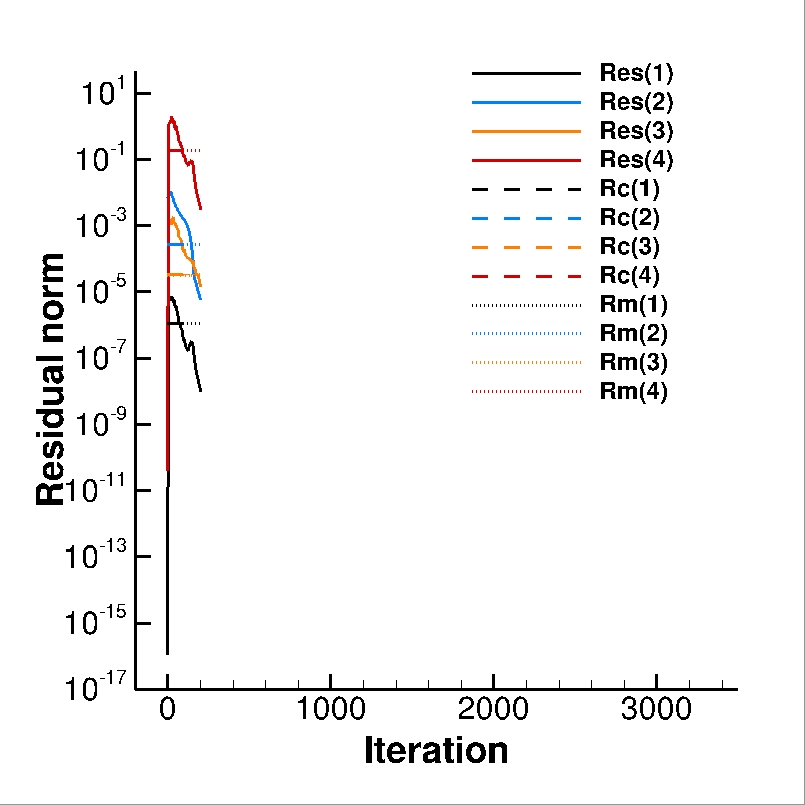}
          \caption{$\epsilon=10^{-6}$.}
          \label{fig:inv_visc_fp_res1_06}
      \end{subfigure}
      \hfill
      \begin{subfigure}[t]{0.32\textwidth}
        \includegraphics[width=\textwidth,trim=2 2 2 2,clip]{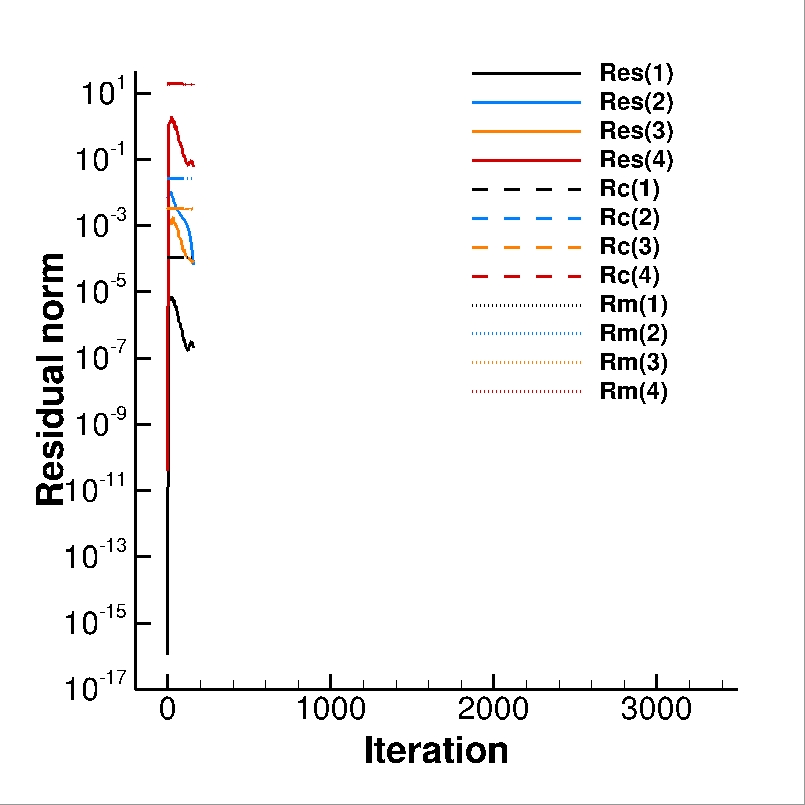}
          \caption{$\epsilon=10^{-4}$.}
          \label{fig:inv_visc_fp_res1_04}
      \end{subfigure}
      \begin{subfigure}[t]{0.32\textwidth}
        \includegraphics[width=\textwidth,trim=2 2 2 2,clip]{./developer_test_res_estimate3/case_visc_fp_tria274x194_4eqns/tol_dw_10m16/myfig_dw1_vs_itr.jpeg}
          \caption{$\epsilon=10^{-16}$.}
          \label{fig:inv_visc_fp_dw1_16}
      \end{subfigure}
      \hfill
      \begin{subfigure}[t]{0.32\textwidth}
        \includegraphics[width=\textwidth,trim=2 2 2 2,clip]{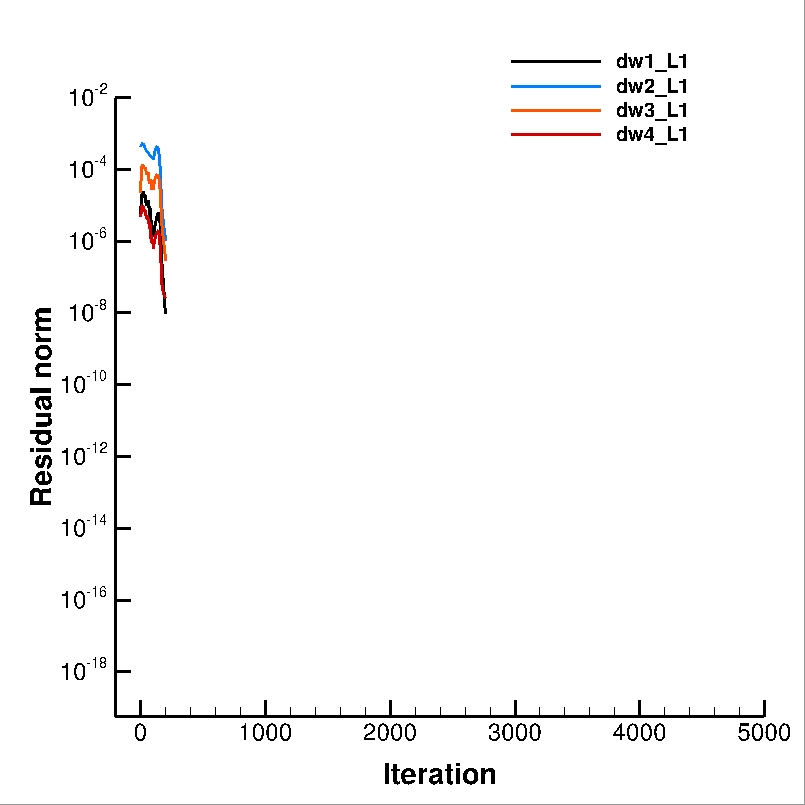}
          \caption{$\epsilon=10^{-6}$.}
          \label{fig:inv_visc_fp_dw1_06}
      \end{subfigure}
      \hfill
      \begin{subfigure}[t]{0.32\textwidth}
        \includegraphics[width=\textwidth,trim=2 2 2 2,clip]{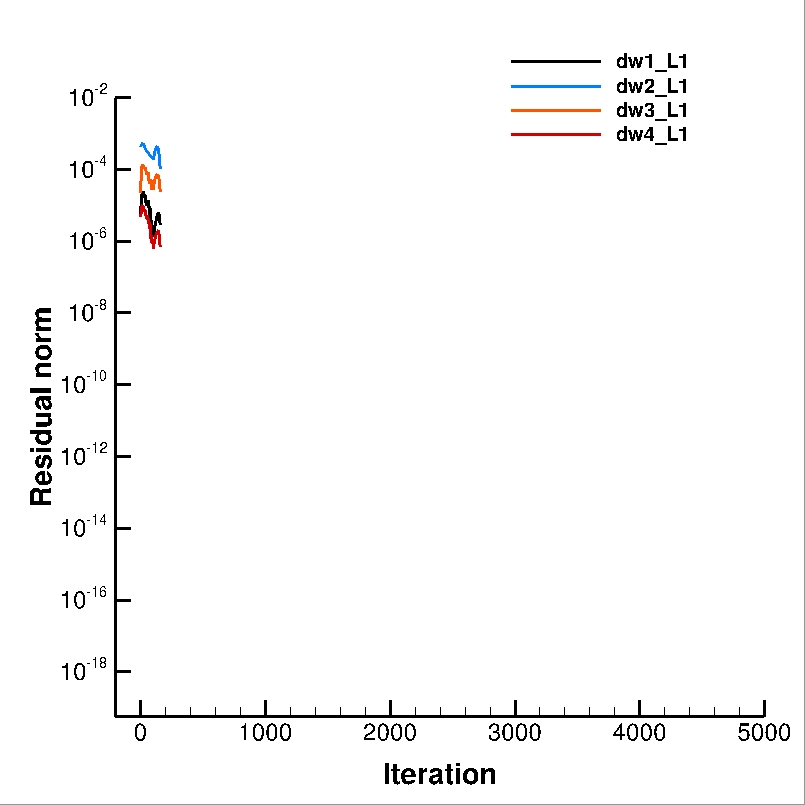}
          \caption{$\epsilon=10^{-4}$.}
          \label{fig:inv_visc_fp_dw1_04}
      \end{subfigure}
      \caption{A flat plate case with different $\epsilon$.}
\label{fig:inv_visc_fp2} 
\end{figure}

\begin{figure}[htbp!]
    \centering
      \begin{subfigure}[t]{0.32\textwidth}
        \includegraphics[width=\textwidth,trim=2 2 2 2,clip]{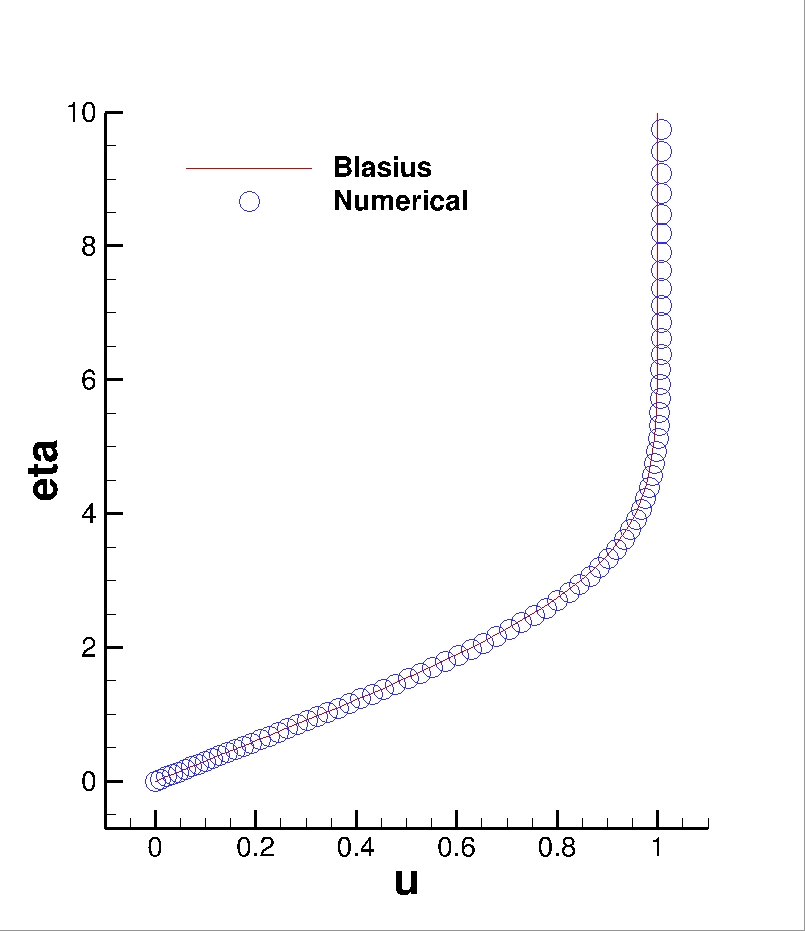}
          \caption{$\epsilon=10^{-16}$.}
          \label{fig:inv_visc_fp_u_16}
      \end{subfigure}
      \hfill
      \begin{subfigure}[t]{0.32\textwidth}
        \includegraphics[width=\textwidth,trim=2 2 2 2,clip]{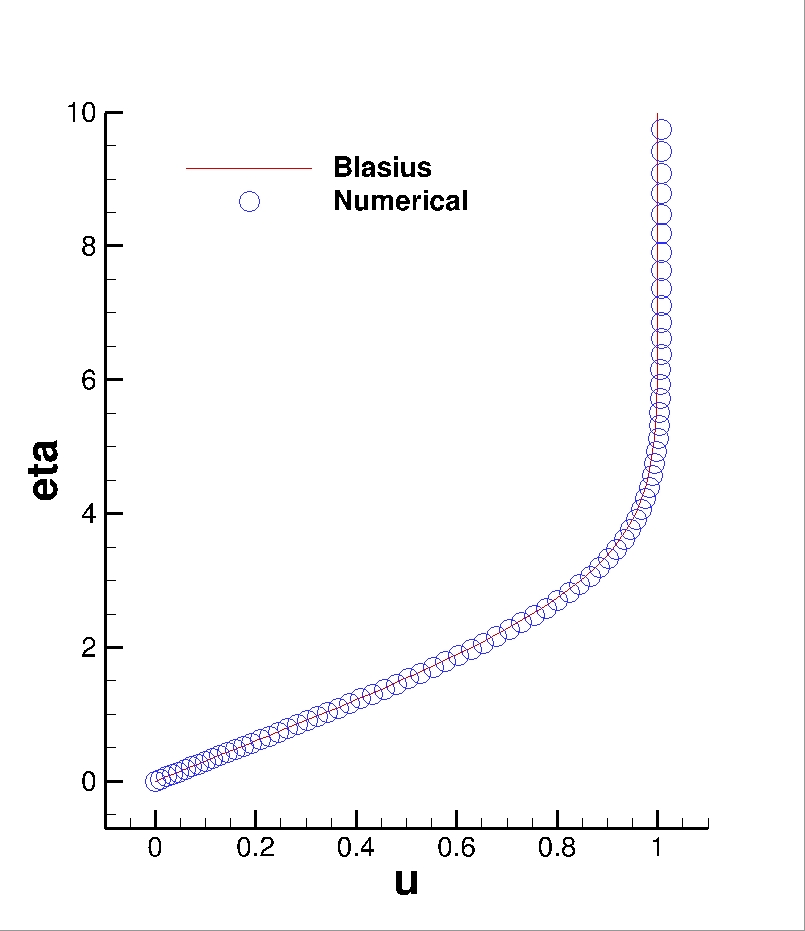}
          \caption{$\epsilon=10^{-6}$.}
          \label{fig:inv_visc_fp_u_06}
      \end{subfigure}
      \hfill
      \begin{subfigure}[t]{0.32\textwidth}
        \includegraphics[width=\textwidth,trim=2 2 2 2,clip]{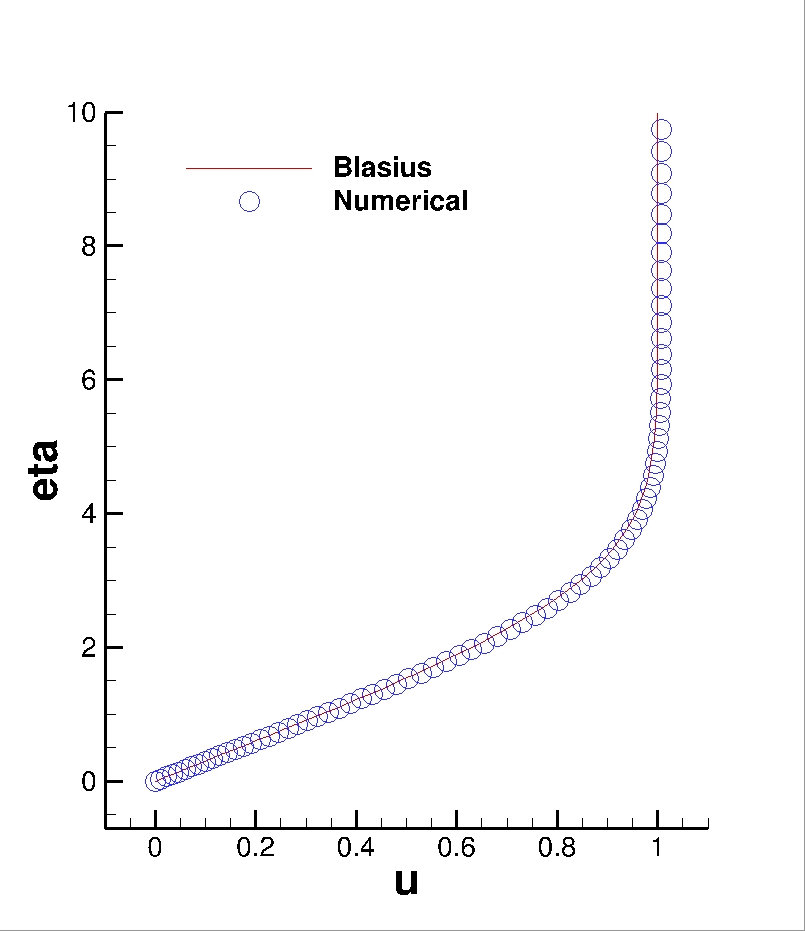}
          \caption{$\epsilon=10^{-4}$.}
          \label{fig:inv_visc_fp_u_04}
      \end{subfigure}
      \begin{subfigure}[t]{0.32\textwidth}
        \includegraphics[width=\textwidth,trim=2 2 2 2,clip]{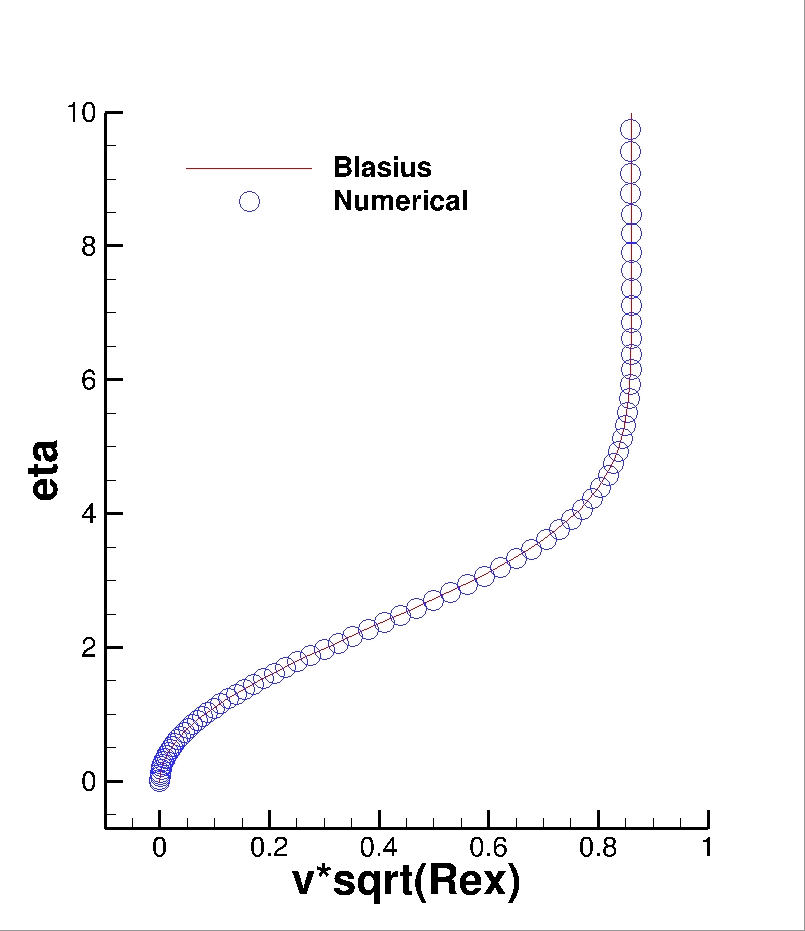}
          \caption{$\epsilon=10^{-16}$.}
          \label{fig:inv_visc_fp_v_16}
      \end{subfigure}
      \hfill
      \begin{subfigure}[t]{0.32\textwidth}
        \includegraphics[width=\textwidth,trim=2 2 2 2,clip]{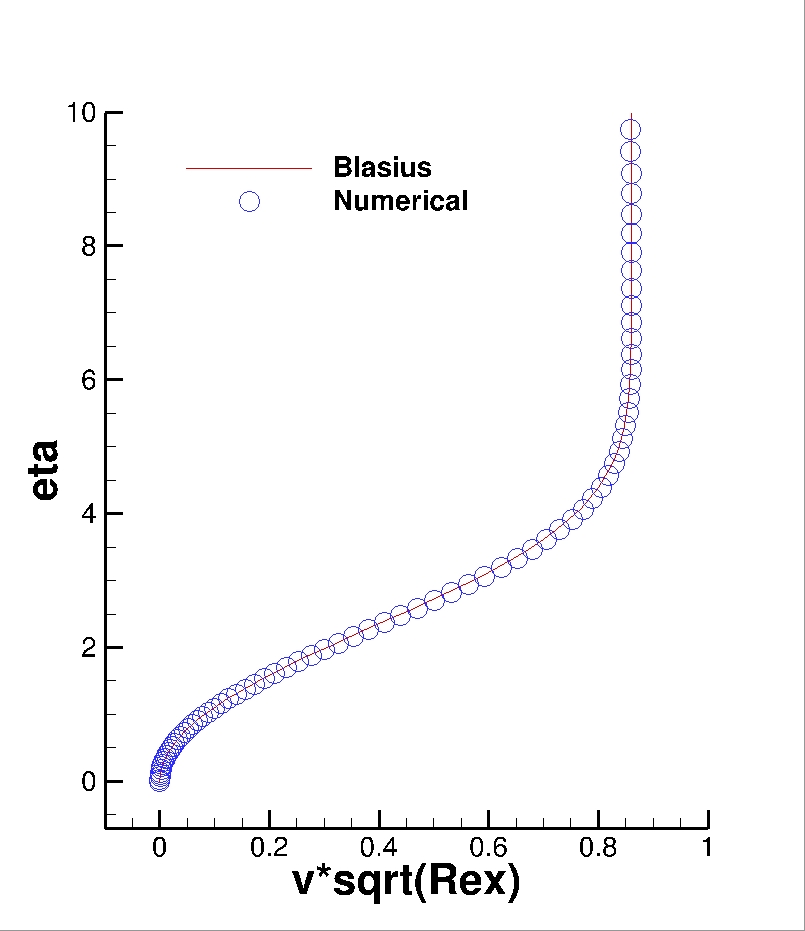}
          \caption{$\epsilon=10^{-6}$.}
          \label{fig:inv_visc_fp_v_06}
      \end{subfigure}
      \hfill
      \begin{subfigure}[t]{0.32\textwidth}
        \includegraphics[width=\textwidth,trim=2 2 2 2,clip]{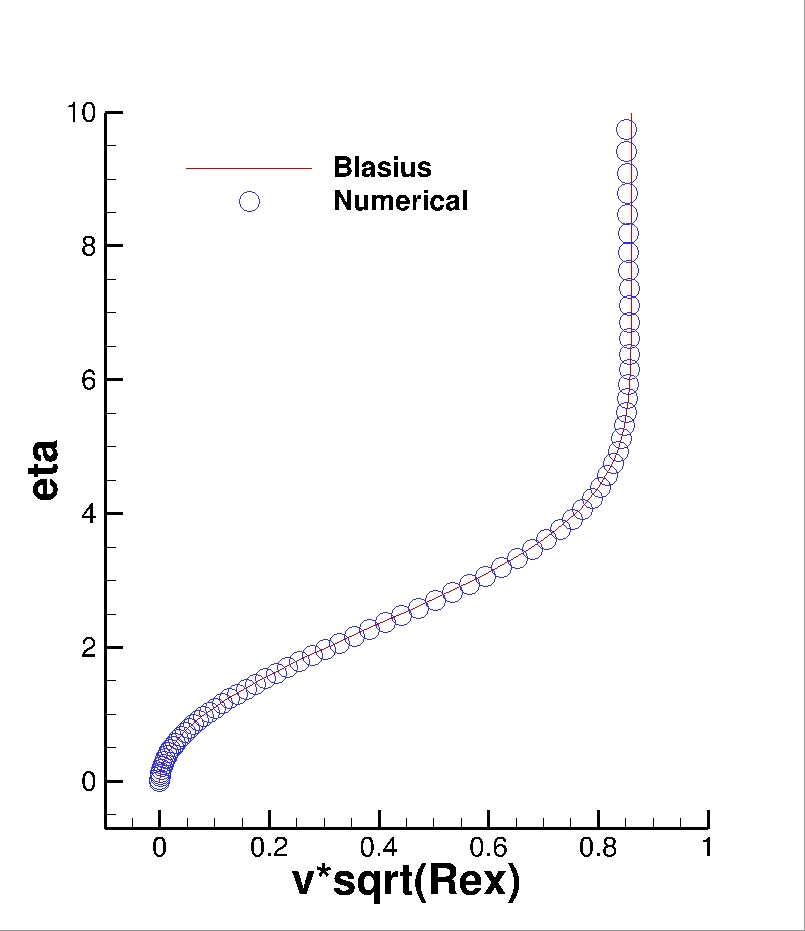}
          \caption{$\epsilon=10^{-4}$.}
          \label{fig:inv_visc_fp_v_04}
      \end{subfigure}
      \begin{subfigure}[t]{0.32\textwidth}
        \includegraphics[width=\textwidth,trim=2 2 2 2,clip]{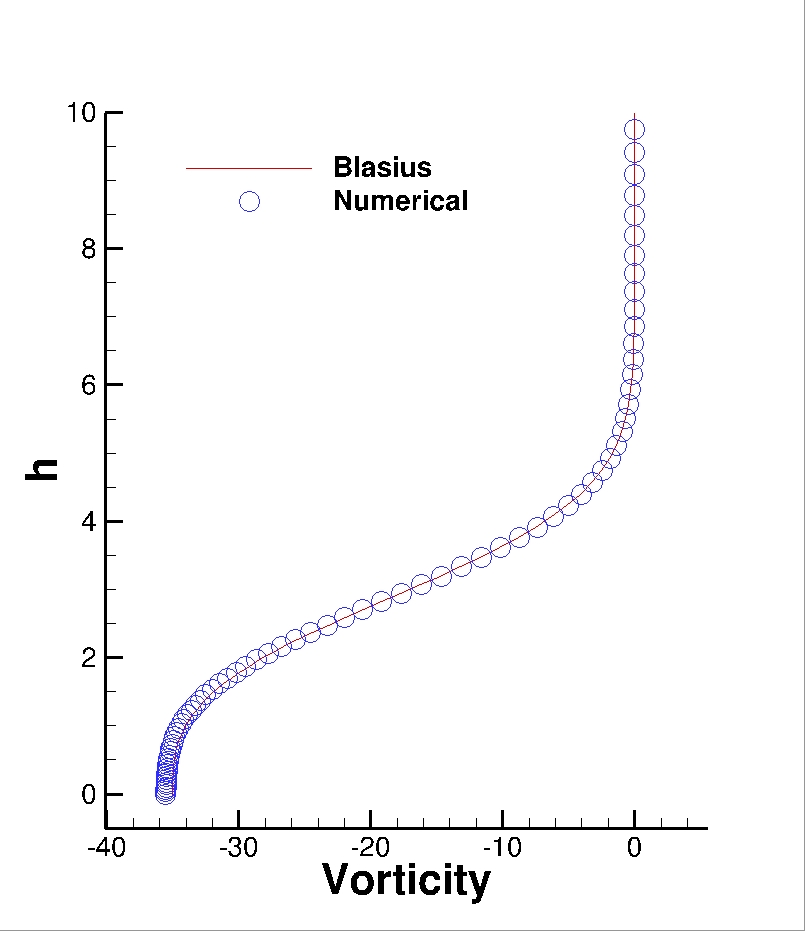}
          \caption{$\epsilon=10^{-16}$.}
          \label{fig:inv_visc_fp_v_16}
      \end{subfigure}
      \hfill
      \begin{subfigure}[t]{0.32\textwidth}
        \includegraphics[width=\textwidth,trim=2 2 2 2,clip]{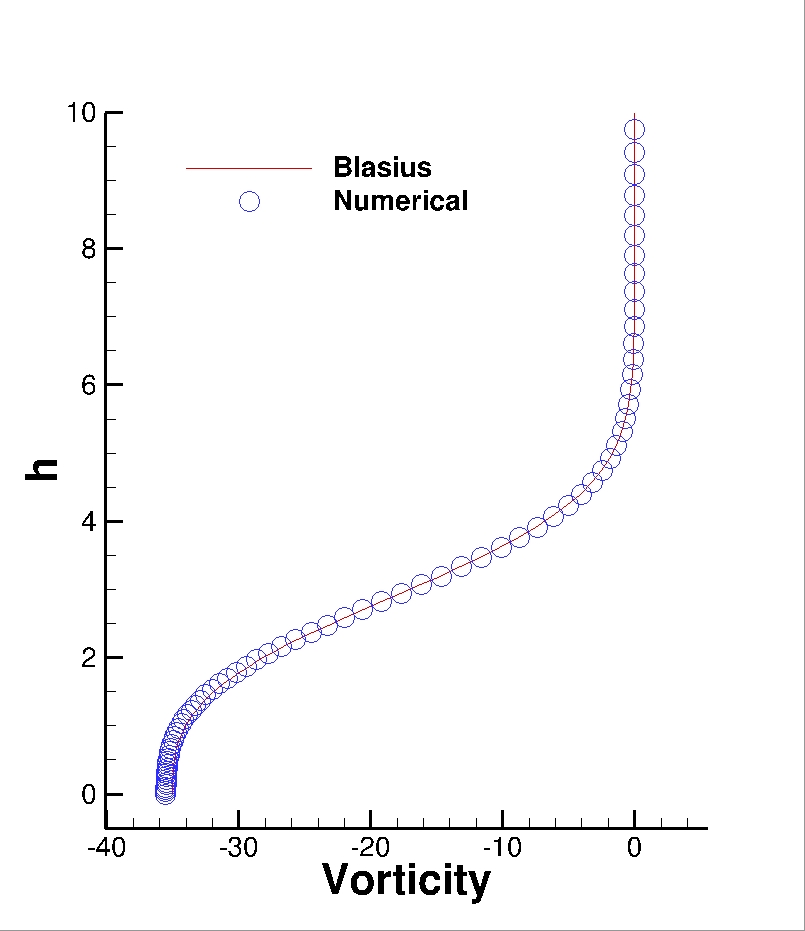}
          \caption{$\epsilon=10^{-6}$.}
          \label{fig:inv_visc_fp_v_06}
      \end{subfigure}
      \hfill
      \begin{subfigure}[t]{0.32\textwidth}
        \includegraphics[width=\textwidth,trim=2 2 2 2,clip]{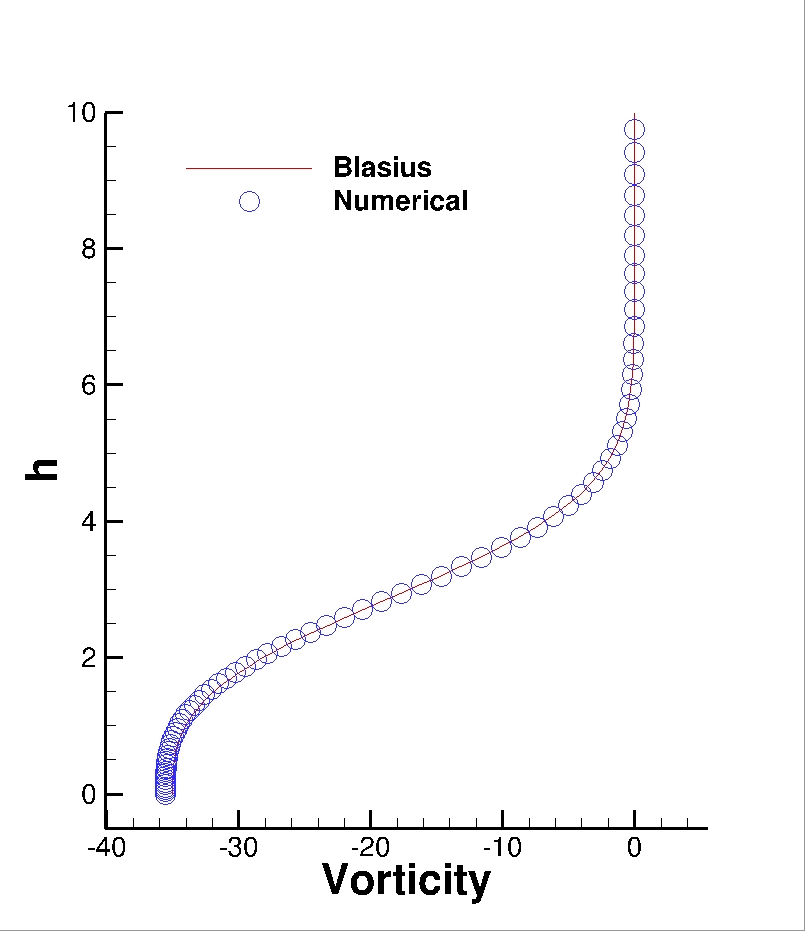}
          \caption{$\epsilon=10^{-4}$.}
          \label{fig:inv_visc_fp_v_04}
      \end{subfigure}
      \caption{A flat plate case with different $\epsilon$. The vertical axis is $\eta$ for all plots; $\eta $ is the normalized coordinate defined by
$\eta = y \sqrt{Re_{x}} / x$, where $Re_{x}$ is the Reynolds number based on the length $x$, i.e., the distance from the leading edge.  }
\label{fig:inv_visc_fp2_sol} 
\end{figure}

For this test case, we explored also the use of larger values of $\epsilon$ to see if these estimates can predict the residual level
for partially converged solutions. Here, the solver is terminated when the iterative solution difference is less than 
a specified value of $\epsilon$. Results are shown in Figure \ref{fig:inv_visc_fp2}. 
As expected, the solver is terminated earlier for a larger value of $\epsilon$. The solutions are very similar as shown in Figure \ref{fig:inv_visc_fp2_sol},
which are sampled along a vertical grid line at $x=0.9$. Therefore, the residuals have been reduced sufficiently for producing accurate solutions. 
However, the residual estimates are not very accurate as can be seen in Figures \ref{fig:inv_visc_fp_res1_06} and \ref{fig:inv_visc_fp_res1_04}.
In fact, we performed the same test for both the airfoil and cylinder cases and found similar results (not shown).
One possible way to use them is to check first the iterative solution difference whether the solution is converged to a desired level and then
use these estimates to confirm that the current residual is small enough (smaller than these estimates); if the residual is larger 
than the estimates, the solver may be stalling.

\section{Conclusions}
\label{conclusions}

We have proposed and demonstrated two techniques to estimate the norm of machine-zero residuals for a given problem.
Numerical results show that these estimates can predict machine-zero residual levels very accurately, at least for the inviscid and viscous flow 
problems considered.
For practical applications, these estimates may be employed to check the convergence of an iterative solver. 
Focusing on how close the magnitude of a current residual is to a machine-zero-residual estimate (instead of how many orders of magnitude it is reduced from an initial residual), we may terminate the iteration when 
the residual norm is, say, five orders of magnitude larger than the machine-zero-residual estimate. 


\addcontentsline{toc}{section}{Acknowledgments}
\section*{Acknowledgments}
 
The second author gratefully acknowledges support from Software CRADLE (part of Hexagon).

\addcontentsline{toc}{section}{References}
\bibliography{../../bibtex_nishikawa_database}

\begin{thebibliography}{1}
\newcommand{\enquote}[1]{``#1''}

\bibitem{liu_nishikawa_aiaa2016-3969}
Liu, Y. and Nishikawa, H., \enquote{Third-Order Inviscid and Second-Order
  Hyperbolic {N}avier-{S}tokes Solvers for Three-Dimensional Inviscid and
  Viscous Flows,} {\em 46th {AIAA} Fluid Dynamics Conference\/}, {AIAA} Paper
  2016-3969, Washington, D.C., 2016.

\bibitem{Roe_JCP_1981}
Roe, P.~L., \enquote{Approximate {R}iemann Solvers, Parameter Vectors, and
  Difference Schemes,} {\em J. Comput. Phys.\/}, Vol.~43, 1981, pp.~357--372.

\bibitem{nishikawa:AIAA2010}
Nishikawa, H., \enquote{Beyond Interface Gradient: A General Principle for
  Constructing Diffusion Schemes,} {\em Proc. of 40th {AIAA} Fluid Dynamics
  Conference and Exhibit\/}, {AIAA} Paper 2010-5093, Chicago, 2010.

\bibitem{WeissMaruszewskiSmith_AIAAJ1999}
Weiss, J.~M., Maruszeski, J.~P., and Smith, W.~A., \enquote{Implicit Solution
  of Preconditioned {N}avier-{S}tokes Equations Using Algebraic Multigrid,}
  {\em {AIAA} J.\/}, Vol.~37, No.~1, 1999, pp.~29--36.

\end{thebibliography}
\bibliographystyle{aiaa}


%

%

%

%
%

\end{document}